\documentstyle[12pt]{article}
\def\mult{\mathop{\rm mult}\nolimits}
\def\C{{\bf C}}
\def\R{{\bf R}}

\def\hx{\widehat {X}}

\def\G{\Gamma}

\def\ha{\widehat {a}}
\def\hb{\widehat {b}}

\def\hf{\widehat {f}}
\def\hg{\widehat {g}}

\def\hu{\widehat {u}}
\def\hv{\widehat {v}}
\def\hx{\widehat {x}}

\def\hq{\widehat {q}}
\def\hz{\widehat {z}}
\def\be{\bar {E}}

\def\p{{\bf P}}
\def\Q{{\bf Q}}

\def\qed{\hfill $\Box$}
\def\be{\begin{equation}}
\def\ee{\end{equation}}

\def\ss{\smallskip}
\def\G{\Gamma}

\def\ul{\underline}
\def\ol{\overline}
\begin{document}

\title{Rational cuspidal plane curves of type $(d,\,d-3)$}

\bigskip

\author{ H. Flenner and  M. Zaidenberg }

\date{}

\maketitle

\bigskip

\begin{abstract} In the previous paper \cite{FlZa 2} we classified
the rational cuspidal plane curves $C$ with a cusp of multiplicity
deg$\,C -2.$  In particular, we showed that
any such curve can be transformed into a line by
Cremona transformations.
Here we do the same for the rational cuspidal
plane curves $C$ with a cusp of multiplicity deg$\,C -3.$
\end{abstract}

\section*{Introduction}

Let $C \subset \p^2$ be a rational cuspidal curve;
that is, it has only irreducible singularities (called
{\it cusps}).
We say that $C$ is of type $(d,\,m)$ if $d =$deg$\,C$
is the degree and
$m = \max_{P \in {\rm Sing}\,C} \{$mult$_P C\}$ is the
maximal multiplicity
of the singular points of $C.$

Topologically, $C$ is a
2-sphere $S^2$ (non-smoothly) embedded into $\p^2.$
Due to the Poincar\'e-Lefschetz dualities, the
complement $X := \p^2 \setminus C$ to $C$ is a
$\Q-$acyclic affine algebraic surface, i.e.
${\widetilde H}_*(X;\,\Q)=0$ (see e.g. \cite{Ra,Fu,Za}). Furthermore,
if $C$ has at least three cusps,  then $X$ is of log-general type, i.e.
${\overline  k}(X)=2,$ where ${\overline k}$ stands for the
logarithmic Kodaira dimension \cite{Wa}.

In \cite{FlZa 1} we conjectured that any $\Q-$acyclic
affine algebraic surface $X$ of log-general type is
rigid in the following sense. Let $V$ be a minimal
smooth projective completion of $X$ by a simple normal
crossing (SNC for short) divisor $D.$  We say that $X$
is {\it rigid} (resp. {\it unobstructed}) if
the pair $(V,\,D)$ has no nontrivial deformations
(resp. if the infinitesimal deformations of the pair $(V,\,D)$ are
unobstructed).

In the particular case when $X = \p^2 \setminus C$
with $C$ as above, the rigidity conjecture would imply
that the curve $C$ itself is projectively
rigid. This means that the only equisingular deformations of
$C$ in $\p^2$ are those
provided by automorphisms of $\p^2;$ in other words,
all of them are projectively equivalent to $C$ (see
\cite[sect.\ 2]{FlZa 2}). In turn,
this would imply that there is only a finite number of
non-equivalent rational cuspidal plane curves of a
given degree with at least three cusps.
Therefore, one may hope to give a classification of such curves.

In \cite{FlZa 2} we obtained a complete list of
rational cuspidal plane curves of type $(d,\,d-2)$
with at least three cusps, and showed that all of them
are projectively rigid and unobstructed.
In the theorem below we do the same for rational
cuspidal plane curves of
type $(d,\,d-3)$ with at least three cusps.

The principal numerical invariant which characterizes a
cusp up to equisingular deformation is its
multiplicity sequence. Recall that, if
$$
V_{n+1}\to V_n\to\dots\to V_1 \to V_0 = \C^2
$$
is a minimal resolution of an irreducible analytic
plane curve germ
$(C,\,{\overline 0})\subset (\C^2,\,{\overline 0}),$
and $(C_i,\,P_i)$ denotes the proper transform of
$(C,\,{\overline 0})$ in $V_i,$ so that $(C_0,\,P_0) =
(C,\,{\overline 0}),$ then $\ul{m}=(m^{(i)})_{i=0}^{n+1},$
where
$m^{(i)} = $mult$_{P_i} C_i,$ is called {\it the
multiplicity sequence} of the germ $(C,\,{\overline
0}).$  Thus, $m^{(i+1)} \le m^{(i)},\,\,\,m^{(n)} \ge 2$
and $m^{(n+1)} =1.$  A multiplicity sequence has the
following characteristic property \cite[(1.2)]{FlZa 2}:

\smallskip

\noindent {\it for any $i=0,\dots,n-1$ either $m^{(i)}
= m^{(i+1)},$ or there exists $k > 0$ such that $i+k
\le n,$ and}
$$m^{(i)} = m^{(i+1)} + \dots + m^{(i+k)} +
m^{(i+k+1)},
\,\,\,\,\,\,\,\,\,\mbox{where}\,\,\,\,\,\,\,\,\,
m^{(i+1)} =\dots = m^{(i+k)}\,.$$
We use the abbreviation $(m_k)$ for a (sub)sequence
$m^{(i+1)} = m^{(i+2)} = \dots =m^{(i+k)}=m.$  Thus, we
present a multiplicity sequence as
$(m^{(1)}_{k_1},\,\dots,m^{(s)}_{k_s})$ with
$m^{(i+1)} < m^{(i)};$ by abuse of notation, we assume
here that $m^{(s)} \ge 2.$
For instance, $(2)$ means an ordinary cusp, and
$(2_3)=(2,\,2,\,2,\,1)$ corresponds to a ramphoid
cusp. With this notation we can formulate our main result as
follows.

\bigskip

\noindent {\bf Theorem}. {\it (a) Let $C \subset
\p^2$ be a rational cuspidal plane curve of type
$(d,\,d-3),\,\,\,d\ge 6,$   with at least three cusps. Then $d = 2k +
3,$ where $k \ge 2,$ and $C$ has exactly three cusps,
of types
$(2k,\,2_k),\,\,(3_k),\,\,(2),$ respectively.

\smallskip

\noindent (b) For each $k \ge 1$ there exists a rational cuspidal plane
curve $C_k$ of degree $d = 2k+3$
with three cusps of types $(2k,\,2_k),$ $(3_k)$ and $(2).$

\smallskip

\noindent (c) Moreover, the curve $C_k$ as in (b)
is unique up to projective equivalence. It can be defined over $\Q.$ }

\bigskip

\noindent {\bf Remarks.} (1) A classification of irreducible plane curves
up to degree $5$ can be found e.g. in \ \cite{Nam}. In particular, there
are, up to projective equivalence, only one rational cuspidal plane quartic
with three cusps ({\it the Steiner quartic}) and only three  rational
cuspidal plane quintic curves with at least three cusps. Two of them have
exactly three cusps, of types $(3),$ $(2_2),$ $(2)$ resp.
$(2_2),$ $(2_2),$ $(2_2),$ and the third one has four cusps of types
$(2_3),$ $(2),$ $(2),$ $(2)$ \ \cite[Thm. 2.3.10]{Nam}.

\smallskip

(2) In his construction of $\Q-$acyclic
surfaces
(see e.g. \cite{tD 1,tD 2}), T. tom Dieck
found certain $(d,\,d-2)-$ and $(d,\,d-3)-$rational
cuspidal curves,
in particular, those listed in the theorem above, as
well as some
other series of rational cuspidal plane curves
(a private communication\footnote{We
are grateful to T. tom Dieck for communicating us the
list of the multiplicity sequences of the constructed
curves.}). Besides a finite number of sporadic examples, 
the curves with at least
three cusps in the list of tom Dieck 
are organized in three  series of $(d,\,d-2)-,$
$(d,\,d-3)-$ and $(d,\,d-4)-$type, respectively. 
It can be checked that all those curves are
rigid and unobstructed. 

Following our methods, T.\ Fenske proved recently that
the only possible numerical data of unobstructed
rational cuspidal plane curves with at least three
cusps and of type $(d,d-4)$ are those from the list of tom Dieck.
He has also classified all
rational cuspidal plane curves of degree 6 \cite{Fe}. 
It turns out that the only
examples with at least 3 cusps are those described in \cite{FlZa 2}.

\smallskip

(3) For a rational cuspidal plane curve $C$ of
type $(d,\,m)$ the inequality $m > d/3$ holds
\cite{MaSa}. Recently, S. Orevkov obtained
a stronger one\footnote{We are grateful to
S. Orevkov for providing us with a preliminary version
of his paper.}: If the complement
$\p^2 \setminus C$ has logarithmic Kodaira dimension 2,
then $d<\alpha m+\beta,$ where $\alpha:=(3+\sqrt{5})/2=2.6180\ldots$
and $\beta :=\alpha-1/\sqrt{5}=2.1708\ldots.$
\smallskip

(4) It was shown in \cite{OrZa 1,OrZa 2}
that a rational cuspidal plane curve with at least ten
cusps cannot be projectively rigid.

\bigskip

Recall the Coolidge--Nagata Problem \cite{Co,Nag}:

\smallskip

\noindent {\it Which rational plane curves can be
transformed into a line by means of Cremona transformations of}
$\p^2?$

\smallskip

\noindent It can be completed by the following question:

\smallskip

\noindent {\it Is this possible, in particular,
for any rational cuspidal plane curve?}

\smallskip

\noindent Under certain restrictions, a positive
answer was given in
\cite{Nag,MKM,MaSa,Ii 2,Ii 3}. It can be verified that the last
question  has a positive answer for the rational cuspidal plane curves
of degree at most five.
In \cite{FlZa 2} we showed that any rational cuspidal
plane curve of type $(d,\,d-2)$ with at least three
cusps is rectifiable. Here we extend this
result to $(d,\,d-3)-$curves. It will turn out to be an immediate
consequence of our construction:

\bigskip
\noindent {\bf Corollary}. {\it Any rational cuspidal
plane curve of type $(d,\,d-3)$ with at least three
cusps is rectifiable, i.e.\ it can be
transformed into a line by means of Cremona
transformations.}

\section{Proofs}

Let $C \subset \p^2$ be a plane curve, and let $V \to
\p^2$ be the minimal embedded resolution of
singularities of $C,$ so that the reduced total
transform $D$ of $C$ in $V$ is an SNC--divisor. By
\cite{FlZa 1}, the
cohomology groups $H^i(\Theta_V\langle \,D\,\rangle)$
of the sheaf of germs of holomorphic vector fields on
$V$ tangent to $D$ control the deformations of the
pair $(V,\,D)$; more precisely,
$ H^0( \Theta_V\langle \, D \, \rangle)$ is the space
of its infinitesimal
automorphisms, $ H^1 ( \Theta_V\langle \, D \,
\rangle)$ is the space of
infinitesimal deformations and $ H^2 ( \Theta_V\langle
\, D \, \rangle)$ gives
the obstructions for extending infinitesimal
deformations.

The surface $X = V \setminus D = \p^2 \setminus C$
being of log-general type, the automorphism group
Aut$X$ is finite \cite{Ii 1}, and hence
$ h^0(\Theta_V\langle \,D\, \rangle) =0.$  Thus, the
holomorphic Euler characteristic of the sheaf
$\Theta_V\langle \, D \, \rangle$ is
$$\chi ( \Theta_V\langle \, D \, \rangle) =
h^2(\Theta_V\langle \,D\,\rangle) -
h^1(\Theta_V\langle \,D\,\rangle).$$
\bigskip

\noindent {\bf Lemma 1.1.} {\it If $C$ is a rational
cuspidal plane curve of type $(d,\,d-3)$ with at least
three cusps, then $h^2(\Theta_V\langle
\,D\,\rangle) = 0,$ that is, $C$ is
unobstructed\footnote{i.e.\ as a plane curve, it has
unobstructed equisingular infinitesimal deformations.}, and so
$\chi=\chi(\Theta_V\langle \,D\,\rangle) \le 0.$ }

\bigskip

\noindent {\it Proof.} Projecting from the cusp of
multiplicity $d-3$
yields a fibration $V \to \p^1,$ which is three--sheeted
when restricted to
the proper transform
of $C.$  Now \cite[(6.3)]{FlZa 1} shows that
$h^2(\Theta_V\langle \, D \,\rangle) = 0.$
Since ${\overline k}(V \setminus D) = 2,$ we also have
$h^0(\Theta_V\langle \, D \,\rangle) = 0.$  Hence
$\chi=-h^1(\Theta_V\langle \, D \,\rangle) \le 0.$
\qed

\bigskip

The next proposition proves part (a) of our main theorem.

\bigskip

\noindent {\bf Proposition 1.1.} {\it The only
possible rational cuspidal plane curves
$C$ of degree $d \ge 6$ with a singular point $Q$ of
multiplicity $d-3$ and at
least three cusps
are those of degree $d = 2k+3,\,\,k=1,\dots,$
with three cusps of types $(2k,\,2_k),$ $(3_k)$ and
$(2).$ Furthermore, these curves are projectively
rigid.}

\bigskip

\noindent {\it Proof.}
By \cite[(2.5)]{FlZa 2} and Lemma 1.1 above, we have:
$$\chi=-3(d-3) + \sum_{P \in {\rm Sing}\,C} \chi_P \le
0\,,\eqno{(R_1)}$$
where
$$\chi_P := \eta_P+\omega_P -1\,,$$ and where, for a
singular point $P \in C$
with the multiplicity sequence
$\ul{m}_P=(m^{(0)},\dots,m^{(k_P)}),$
$$
\eta_P = \sum\limits_{i=0}^{k_P}
(m^{(i)}-1)\,\qquad{\rm
and}
\qquad
\omega_P = \sum\limits_{i=1}^{k_P} (\lceil{m^{(i-1)}
\over m^{(i)}} \rceil - 1)\,
$$
(for $a \in \R,\,\,\lceil {a} \rceil$ denotes the
smallest integer $\ge a$).

Observe that, by the Bezout theorem, $m_P^{(0)} +
m_P^{(1)} \le d$ and
$m_P^{(0)} + m_Q^{(0)} \le d.$  Thus
$$
\mbox{for} \quad P\neq Q
\quad\mbox{we have}\quad
 m_P^{(0)} \le 3;\quad
\mbox{moreover we have} \quad m_Q^{(1)} \le 2,
$$
since otherwise the tangent line $T_Q C$
would have the only point $Q$ in common with $C,$ and
so, $C \setminus T_Q C$
would be an affine rational cuspidal plane curve with one
point at infinity and with two cusps.
But by the Lin-Zaidenberg Theorem \cite{LiZa},
up to biregular automorphisms
of the affine plane $\C^2,$ the
only irreducible simply connected affine plane
curves are the curves
$\Gamma_{k,\,l} = \{x^k - y^l = 0\},$ where $1 \le k
\le l,$ and $(k,\,l) = 1.$  Hence, such a curve cannot have two
cusps. Using the above restriction and the characteristic
property of a multiplicity sequence cited above we obtain the
following possibilities for the multiplicity sequence $\ul{m}_P$ at
a singular point $P$:
$$
\begin{array}{l}%
\ul{m}_Q=(d-3)\mbox{ or }(d-3,2),\\[2pt]
\ul{m}_P=(2_a) \mbox{ or } (3_a)\mbox{ or }(3_a,2)\quad
\mbox{for } P\ne Q .
\end{array}\eqno{(R_2)}
$$
For different possible types of cusps of $C$ we have:

\medskip

\noindent (a) If $Q \in $\,Sing$\,C$ has the
multiplicity sequence $(d-3),$
then
$$
\eta_Q = d-4,\quad\omega_Q = d-4\quad\mbox{
and so}\quad\chi_Q = 2d-9.
$$

\medskip

\noindent (b) If $Q \in $\,Sing$\,C$ has the
multiplicity sequence $(d-3,\,2_a)$ then, by the same
characteristic property \cite[(1.2)]{FlZa 2},
$$
\quad \mbox{either}\quad d-3\le 2a \quad\mbox{is even or}
\quad d-3 = 2a+1.\leqno (*)
$$
In any case
$$
\eta_Q = d-4 + a,\quad\omega_Q =
\lceil{d-3\over 2} \rceil\quad\mbox{and so}\quad\chi_Q = d-5 + a +
\lceil{d-3\over 2} \rceil.
$$

\medskip

\noindent (c) If $P \in $\,Sing$\,C$ has the
multiplicity sequence $(2_a),$
then
$$
\eta_P = a,\quad\omega_P =
1\quad
\mbox{and so}\quad\chi_P = a.
$$

\medskip

\noindent (d) If $P \in $\,Sing$\,C$ has the
multiplicity sequence $(3_a),$
then
$$
\eta_P = 2a,\quad\omega_P = 2 \quad\mbox{and so}\quad\chi_P =2a+1.
$$

\medskip

\noindent (e) If $P \in $\,Sing$\,C$ has the
multiplicity sequence $(3_a,\,2),$ then
$$
\eta_P =2a+1,\quad\omega_P =2\quad\mbox{and so}\quad\chi_P =2a+2.
$$

\medskip

Furthermore, since $C$ rational, by the genus formula,
we have
$$
{d-1 \choose 2} = \sum_{P \in {\rm Sing}\,C}\delta_P
\quad\mbox{where}\quad
\delta_P:= \sum\limits_{i=1}^{k_P} {m_P^{(i)}\choose 2}\,.
$$
Since $m_Q^{(0)} = d-3,$ we get
$$
{d-1 \choose 2} = {d-3 \choose 2} + \sum_{(P,\,i)
\neq (Q,\,0)} {m_P^{(i)} \choose 2}\,,
$$
or, equivalently,
$$
2d-5 = \sum_{(P,\,i) \neq (Q,\,0)}
{m_P^{(i)}(m_P^{(i)}-1)\over 2}\,.\eqno{(R_3)}
$$

At last, consider the projection $\pi_Q\,:\,C \to
\p^1$ from the point $Q.$  By the Riemann-Hurwitz
Formula, it has at most four branching points.
This gives the restriction (see \cite[(3.1)]{FlZa 2})
$$m^{(1)}_Q - 1 + \sum\limits_{P\neq Q} (m^{(0)}_P -
1) \le 4\,. \eqno{(R_4)}$$
Thus, if the curve $C$ has the numerical data
$$
[(d-3,\,2_{a_1}),\,(2_{a_2}),
\dots,(2_{a_k}),(3_{b_1}),\dots,(3_{b_l}),
(3_{c_1},\,2),\dots,(3_{c_m},\,2)]\,,
$$
then $k+2(l+m) \le 4.$  Hence, either $l+m=0$ and $3\le
k \le 4,$ or
$l+m=1$ and $k= 2,$ or $l+m=2$ and $k =0.$

Taking into account the above restrictions $(R_2) -
(R_4)$ and $(*)$ from (b), the list of all possible data of
rational cuspidal plane curves
$C$ of degree $d\ge 6$ with a point of multiplicity
$d-3$ and at least $3$ cusps is as follows, where $a,\,b,\,c,\,e > 0$:

\be
[(d-3),\,\,\,(2_a),\,\,\,(2_b)]
\quad \mbox{where} \quad a+b = 2d-5
\ee
\be [(d-3),\,\,\,(2_a),\,\,\,(3_b)]
\quad \mbox{where} \quad a+3b = 2d-5 \ee
\be [(d-3),\,\,\,(3_a),\,\,\,(3_b)]
\quad \mbox{where} \quad 3a+3b = 2d-5 \ee
\be [(d-3),\,\,\,(3_a,\,2),\,\,\,(2_b)]
 \quad \mbox{where} \quad 3a+b = 2d-6 \ee
\be [(d-3),\,\,\,(3_a,\,2),\,\,\,(3_b)]
 \quad \mbox{where} \quad 3a+3b = 2d-6 \ee
\be [(d-3),\,\,\,(3_a,\,2),\,\,\,(3_b,\,2)]
 \quad \mbox{where} \quad 3a+3b = 2d-7 \ee
\be [(d-3,\,2_a),\,\,\,(2_b),\,\,\,(2_c)]
 \quad \mbox{where} \quad a+b+c =
2d-5\mbox{ and }  (*)\mbox{ holds}\ee
\be [(d-3,\,2_a),\,\,\,(3_b),\,\,\,(2_c)]
 \quad \mbox{where} \quad a+3b +c=
2d-5\mbox{ and }  (*)\mbox{ holds} \ee
\be [(d-3,\,2_a),\,\,\,(3_b,\,2),\,\,\,(2_c)]
 \quad \mbox{where} \quad a+3b +c=
2d-6\mbox{ and }  (*)\mbox{ holds}\ee
\be [(d-3),\,\,\,(2_a),\,\,\,(2_b),\,\,\,(2_c)]
 \quad \mbox{where} \quad
a+b +c= 2d-5\ee
\be [(d-3),\,\,\,(3_a),\,\,\,(2_b),\,\,\,(2_c)]
 \quad \mbox{where} \quad 3a+b +c= 2d-5\ee
\be [(d-3),\,\,\,(3_a,\,2),\,\,\,(2_b),\,\,\,(2_c)]
 \quad \mbox{where} \quad 3a+b +c= 2d-6\ee
\be [(d-3,\,2_a),\,\,\,(2_b),\,\,\,(2_c),\,\,\,(2_e)]
 \quad \mbox{where}\quad a+b +c+e=
2d-5\mbox{ and }  (*)\mbox{ holds}\ee
\be
[(d-3),\,\,\,(2_a),\,\,\,(2_b),\,\,\,(2_c),\,\,\,(2_e)]
 \quad \mbox{where} \quad a+b +c+e= 2d-5.\ee

We will examine case by case, computing $\chi=\chi (
\Theta_V\langle \, D \, \rangle).$  The genus formula
and the restriction $(R_1)$ $\chi\le 0$
provided by Lemma 1.1 will allow
to eliminate all the cases but one, namely, a subcase
of (8).

\medskip

\noindent Case (1):
$[(d-3),\,\,\,(2_a),\,\,\,(2_b)]$  where $a+b =
2d-5.$
By $(R_1),$ we have
$\chi = (-3d+9)+(2d-9) + (a + b) = d-5 \le 0,$ a
contradiction.

\smallskip

\noindent Case (2): $[(d-3),\,\,\,(3_b),\,\,\,(2_a)]$
where $a+3b = 2d-5.$
We have
$\chi = (-3d+9)+(2d-9) + (a + 2b +1) = d-4-b \le 0,$
i.e.\ $b \ge d-4.$
On the other hand, $2d -5= a + 3b  \ge 3b + 1,$ whence
$b \le {2\over 3}d - 2.$  Therefore, $d - 4 \le {2\over
3}d - 2,$ i.e.\ $d \le 6.$  In the case $d = 6$ the only
possibility would be $[(3),\,\,(2),\,\,(3_2)].$
Projecting from the cusp with the  multiplicity
sequence $(3_2),$ we get a
contradiction to the Hurwitz formula (see ($R_4$)).

\smallskip

\noindent Case (3):
$[(d-3),\,\,\,(3_a),\,\,\,(3_b)]$  where $ 3a+3b =
2d-5.$
We have
$\chi = (-3d+9)+(2d-9) + (2a +1 + 2b + 1) = {d-4 \over
3} \le 0,$ i.e.\
$d \le 4,$ and we are done.

\smallskip

\noindent Case (4):
$[(d-3),\,\,\,(3_a,\,2),\,\,\,(2_b)]$ where $3a+b =
2d-6.$
We have
$\chi = (-3d+9)+(2d-9) + (2a +2 + b) = d - 4 - a \le
0,$ i.e.\ $a \ge d-4.$
But $2d-6 =  3a + b \ge 3a + 1,$ whence $a \le {2\over
3}d - {7\over 3},$
and thus $d - 4 \le {2\over 3}d - {7\over 3},$ or $d
\le 5,$ a contradiction.

\smallskip

\noindent Case (5):
$[(d-3),\,\,\,(3_a,\,2),\,\,\,(3_b)]$
where $3a+3b = 2d-6.$
We have
$\chi = (-3d+9)+(2d-9) + (2a +2 + 2b+1) = {d\over 3} -
1 \le 0,$ i.e.\ $d\le 3,$
which is impossible.

\smallskip

\noindent Case (6):
$[(d-3),\,\,\,(3_a,\,2),\,\,\,(3_b,\,2)]$ where
$3a+3b = 2d-7.$  We have
$\chi = (-3d+9)+(2d-9) + (2a + 2b+4)= {d\over 3} -
{2\over 3} \le 0,$ which is impossible.

\smallskip

\noindent Case (7):
$[(d-3,\,2_a),\,\,\,(2_b),\,\,\,(2_c)]$ where $a+b+c
= 2d-5$ and $(*)$ holds. We have
$\chi = (-3d+9)+ (d-5 +a +\lceil{d-3
\over 2} \rceil) + (b+c) = \lceil{d-3
\over 2} \rceil - 1  \le 0,$ or $d \le 5,$ and we are
done.

\smallskip

\noindent Case (8):
$[(d-3,\,2_a),\,\,\,(3_b),\,\,\,(2_c)]$ where $a+3b
+c= 2d-5$ and $(*)$ holds. We have
$\chi = (-3d+9)+ (d-5 +a +\lceil{d-3\over 2} \rceil) +
(2b+1 +c) = \lceil{d-3
\over 2} \rceil - b  \le 0,$ i.e.\ $b \ge \lceil{d-3
\over 2} \rceil.$

If $d-3$ is odd, then we get $2d-5 = a + 3b + c \ge 3b
+ 1 + {d-4 \over 2},$
as $a = {d-4 \over 2}$ by $(*).$  Hence, $b \le {d\over
2} - {4\over 3}.$
This leads to $ \lceil{d-3
\over 2} \rceil = {d-2 \over 2}\le {d\over 2} -
{4\over3},$ which is a contradiction.

If $d-3$ is even, then by $(*)$ we get $2d-5= a + 3b
+ c \ge 3b + 1 + {d-3 \over 2},$ hence
$b \le {d\over 2} - {3\over 2}.$  Thus, $\lceil{d-3
\over 2} \rceil =
{d-3 \over 2} \le b \le {d - 3 \over 2},$ which is only
possible if
$c = 1,\,\,a=b={d-3 \over 2}.$  With $k:= {d - 3
\over 2}$ we obtain that
$d = 2k+3,\,\,a=b=k$ and $c = 1;$ that is, $C$ is as
in the proposition.
Observe that in this case $\chi = 0,$ and so
$h^1(\Theta_V\langle \, D \,\rangle) = 0.$  Together with Lemma 1.1
this proves that the corresponding curve $C$ is projectively rigid
and unobstructed (see [FZ 2, Sect. 2]).

\smallskip

\noindent Case (9):
$[(d-3,\,2_a),\,\,\,(3_b,\,2),\,\,\,(2_c)]$ where
$a+3b +c= 2d-6$ and $(*)$ holds. We have
$\chi = (-3d+9)+ (d-5 +a +\lceil{d-3\over 2} \rceil) +
(2b+2 +c) =
\lceil{d-3 \over 2} \rceil - b  \le 0,$ which gives $b
\ge \lceil{d-3
\over 2} \rceil.$

If $d-3$ is odd, then we get $2d-6 = a + 3b + c \ge 3b
+ 1 + {d-4 \over 2},$
as $a = {d-4 \over 2}$ by $(*).$  Thus, $b \le {d\over
2} - {5\over 3},$ and so
we have ${d-2\over 2}\le {d\over 2} - {5\over 3},$
which is a contradiction.

If $d-3$ is even, then we get $2d-6 = a + 3b + c  \ge
3b + 1 + {d-3\over 2}.$  Hence, $b \le {d\over 2} -
{11\over 6}.$  This yields
${d-3\over 2}\le {d\over 2} - {11\over 6},$ which
again gives a contradiction.

\smallskip

\noindent Case (10):
$[(d-3),\,\,\,(2_a),\,\,\,(2_b),\,\,\,(2_c)]$ where
$a+b +c= 2d-5.$  We have $\chi = (-3d+9)+(2d-9) +(a +b
+c) = d-5 \le 0,$ and we are done.

\smallskip

\noindent Case (11) resp. (12), (13), (14) can be
ruled out by the same computations as in case (2)
resp. (4), (7), (10).
This completes the proof of Proposition 1.1.
\qed

\bigskip

For the proof of part (b) and (c) the main theorem we need the
following facts.

\bigskip

\noindent {\bf Lemma 1.2.} {\it
Let $(C,\,{\ol 0}),\ (D,\,{\ol 0})\subseteq(\C^2,\,{\ol 0})$ be two curve
singularities which
have no component in common. Then the following hold.

\smallskip

\noindent (a) $(CD)_{\ol 0}=\sum_P\mult_PC \mult_PD,$ where the sum is
taken over ${\ol 0}$ and all its infinitesimally near points.

\smallskip

\noindent (b) Assume that $(D,\,{\ol 0})$ is a smooth germ and $(C,\,{\ol 0}=
)$
is a cusp with
the multiplicity sequence $\ul{m}=(m^{(0)},\ldots, m^{(n)}).$ Then
$(CD)_{\ol 0}=
m^{(0)}+\ldots+ m^{(s)}$ for some $s\ge 0,$ where $m^{(0)}=\ldots=
m^{(s-1)}.$

\smallskip

\noindent (c) Let $\pi:X\to \C^2$ be the blow up at ${\ol 0}.$  Denote by
$E\subseteq X$ the exceptional curve, and by $C'$ the proper
transform of $C.$ Then}
$$
\mult_{\ol 0}C=\sum_{P\in E}(EC')_P.
$$
\smallskip

\noindent  {\it Proof.} The statements (a) and (c) are well known
(see e.g.\ [Co]), whereas (b) is shown in
[FlZa 2, (1.4)]. \qed

\bigskip

The next result proves part (b) and (c) of the main theorem as well as
the corollary from the introduction.

\medskip

\noindent {\bf Proposition 1.2.} \it (a) For each $k \ge
1$ there exists a rational cuspidal plane curve $C_k$ of degree $d = 2k+3$
with three cusps $Q_k,\ P_k,\ R_k$ of types $(2k,\,2_k),$ $(3_k)$ and
$(2),$ respectively.

(b) $C_k$ is unique up to a projective
transformation of the plane.

(c) $C_k$ is defined over $\Q$.

(d) $C_k$ is rectifiable.
\bigskip

\rm
\noindent {\it Proof.} We proceed by induction on $k.$
Namely, given a curve $C_k$ as in (a), we find a Cremona transformation
$\psi_k\,:\,\p^2\to\p^2$ such that the proper transform
$C_{k+1}=\psi_k(C_k)$ of $C_k$ under $\psi_k$ is a cuspidal curve of
degree
$2k+5$ with three cusps of type
$(2k+2,\,2_{k+1}),\,\, (3_{k+1}),\,\, (2).$  Hence the existence
follows. This construction will also show that (b)--(d) hold.

We start with the rational cuspidal cubic $C_0\subseteq\p^2$
given by the equation $x^2z=y^3.$  Observe that $C_0$ is rectifiable.
It has a simple cusp at $R_0:=(0:0:1)$ and the only inflectional
tangent line $\ell_0$ at $P_0:=(1:0:0);$ that is,
$\ell_0\cdot C_0=3P_0.$
Fix an arbitrary point\footnote{Observe that the projective
transformation group
$(x : y : z) \longmapsto (t^3x : t^2y : t^6z),\,\,\,t \in \C^*,$ acts
transitively in $C_0 \setminus \{P_0,\,R_0\}.$ } $Q_0\in C_0 \setminus
\{P_0,\,R_0\}.$  Let $t_0$ be the tangent line to $C_0$ at $Q_0;$ then
we have $t_0\cdot C_0=2Q_0 + S_0,$ where, as it is easily seen, $S_0 \in C=
_0$
is different from $P_0,\,Q_0$ and $R_0.$  Let $Q_0^*$ denote the
intersection point
$l_0 \cap t_0;$ clearly, $Q_0^* \notin C_0.$

Let, for a given $k > 0,$ $C_k$ denotes a curve with the cusps
$Q_k,\ P_k,\ R_k$ as in the proposition, and let
$C_0$ be the rational cubic with the distinguished points
$Q_0,\ P_0,\ R_0, \ S_0$ as described above. For $k>0$ let
$t_k$ be the tangent line of $C_k$ at $Q_k,$ and $\ell_k$ be the
line $\ol{P_kQ_k},$ whereas for $k=0$ we choose $t_0$
and $\ell_0$ as above. In any case, using Bezout's Theorem and Lemma
1.2, we have
$$
\ell_k\cdot C_k=(d-3)Q_k + 3P_k,\quad\mbox{and}\quad
t_k\cdot C_k=(d-1)Q_k + S_k\,,
$$
where $S_k \in C_k$ is different from $P_k,\,Q_k$ and $R_k.$
Indeed, the line $t_k$ intersects $C_k$ at the point $Q_k$ with multiplicity
$d-1$ if $k>1$ (see Lemma 1.2 (b)) or $k=0.$  To show that this is
also true for $k=1,$
assume that $t_1$ and $C_1$ only intersect in $Q_1$ with
$(t_1C_1)_{Q_1}=d=5.$  The linear projection from $Q_1$
yields a 3-sheeted covering of the normalization of
$C_1$ onto $\p^1.$  By the Riemann-Hurwitz formula, it must have four
ramification points. But since $(t_1C_1)_{Q_1}=d=5,$ the point $Q_1$
would be a ramification point of index $\ge 2$ (see Lemma 1.2(a)),
and so we would have three
ramification points $Q_1,\, P_1,\, R_1$ of indices
$2,\, 2,\, 1,$ respectively, which is a contradiction.

Hence, for any $k\ge 0$ there is exactly one
further intersection point $S_k\in C_k\cap t_k$ with
$(t_kC_k)_{S_k}=1.$

Let $\sigma_k\,:\,X_k\to\p^2$ be the blow up at the point
$t_k\cap\ell_k,$ which is $Q_k$ for $k>0$ and $Q_k^*$ for $k=0.$
Denote by $C'_k,$ $\ell_k',$
$t_k'$ the proper transforms in $X_k$ of the curves $C_k,$ $\ell_k,$ $t_k,$
respectively. Then $X_k\simeq \Sigma_1$ is a Hirzebruch surface with a
ruling $\pi_k\,:\,X_k \to \p^1$ given by
the pencil of lines through $Q_k$ resp.\ $Q_0^*,$ and with the exceptional
section $E_k = \sigma_k^{-1}(Q_k),\,k>0,$ resp.
$E_0 = \sigma_0^{-1}(Q_0^*),$ where $E_k^2 = -1.$  Thus, $\ell_k',$ $t_k=
'$
are fibres of this ruling. By construction, the restriction
$\pi_k\,|\,C_k'\,:\,C_k' \to \p^1$ is 3-sheeted, and we have
$$
\ell_k'\cdot C_k'=3P_k',\quad t_k'\cdot C_k'=2Q_k' + S_k',
\quad\mbox{and}\quad E_k'\cdot C_k'=(d-3)Q_k' = 2kQ_k'\,,
$$
where $P_k',\,Q_k',\,R_k'$ and $S_k'$ are the points of $C_k'$
infinitesimally near to $P_k,\,Q_k,\,R_k$ and $S_k \in C_k,$ respectively
(indeed, by Lemma 1.2(c), we have $(E_k'C_k')_{Q_k'} =
$mult$_{Q_k^*}C_k = d-3,$ where for $k > 0$ we set $Q_k^* = Q_k$).
Clearly,
for $k > 0$ $P_k',\,Q_k'$ and $R_k'$ are cusps of $C_k'$ of types
$(3_k),\,\,(2_k)$ and $(2),$ respectively, whereas $S_k'$ is a smooth point.

Next we perform two elementary transformations\footnote{Recall that an
elementary transformation of a ruled surface consists in blowing up at a
point of a given irreducible fibre followed by the contraction of the
proper
transform of this fibre.} of $X_k,$ one at the
point $S_k'$ and the other one at the intersection point $T_k':=\{E_k\cap
\ell_k'\}.$  We arrive at a new Hirzebruch surface $X_{k+1}\simeq\Sigma_1,$
with the exceptional section
$E_{k+1}$ being the proper transform of $E_k$ (indeed, since we perform
elementary transformations at the points $S_k \notin E_k$ and $T_k' \in E_k,=
$
we have $E_{k+1}^2=E_k^2=-1$). Denote by $C_{k+1}'$ the
proper transform of
$C_k',$ and by $t_{k+1}',$ $\ell_{k+1}'$ the fibres of the ruling
$\pi_{k+1}\,:\,X_{k+1}\to \p^1$ which replace $t_k'$ resp.\ $\ell_k'.$
Using formal properties
of the blowing up/down process we obtain, once again, the relations
$$
\ell_{k+1}'\cdot C_{k+1}'=3P_{k+1}',\quad t_{k+1}'\cdot C_{k+1}'=2Q_{k+1=
}'
+ S_{k+1}',
\quad\mbox{and}\quad E_{k+1}'\cdot C_{k+1}'= 2(k+1)Q_{k+1}'\,,
$$
where $P_{k+1}',\,Q_{k+1}',\,R_{k+1}'$ and $S_{k+1}'$ are the points of
$C_{k+1}'$ infinitesimally near to $P_k',\,Q_k',\,R_k'$ and $S_k' \in
C_k',$ respectively. It is easily seen that $P_{k+1}'$ resp.
$Q_{k+1}',\,\,\,R_{k+1}'$
are cusps of $C_{k+1}'$ of types
$(3_{k+1}),\,\,(2_{k+1})$ and $(2),$ respectively, whereas $S_{k+1}'$ is a
smooth point.

Blowing down the exceptional curve $E_{k+1}' \subset X_{k+1}$
we arrive again at $\p^2.$ Denote the images of $C_{k+1}',.$$  Q_{k+1}',$
$P_{k+1}',.$$  R_{k+1}'$ resp. by $C_{k+1},\ Q_{k+1},$
$P_{k+1},$ $R_{k+1}.$  We have constructed a rational cuspidal plane curve
$C_{k+1}$ which has
cusps at $Q_{k+1},\ P_{k+1},\ R_{k+1}$ with multiplicity sequences
$(2(k+1),\, 2_{k+1}),\,\, (3_{k+1}),\,\, (2),$ respectively (see Lemma
1.2(c)). This completes the proof of existence.

Note that the birational transformation
$\psi_k:\p^2\to\p^2$, by which we obtained $C_{k+1}=\psi_k(C_k)$
from $C_k$, is just the Cremona transformation in the points $S_k$,
$Q_k$ and the intersection point $E_k\cap \ell_k'$, which is
infinitesimally near to
$Q_k$. This transformation only depends upon $Q_k$, $S_k$ and the line
$\ell_k;$ we denote it by
$\psi(S_k,Q_k,\ell_k):=\psi_k.$ The inverse $\psi_k^{-1}$ is the
transformation $\varphi_k=\psi(P_{k+1}, Q_{k+1},t_{k+1})$. Therefore, the
curve $C_k$ is always transformable into the cuspidal cubic, and thus also
into a line, by means of Cremona transformations, proving (d).  In order
to show (c) we note that, moreover, so constructed $C_k,$ as well
as $P_k$, $Q_k$, $R_k$ and $S_k,$ are defined over $\Q,$ as follows by an
easy induction.

\bigskip

Finally, let us show that the curve $C_k$ is uniquely determined up
to a projective transformation of the plane. We will again proceed by
induction on $k$. Clearly, the cuspidal cubic is uniquely determined up
to a projective transformation. Assume that uniqueness is shown for the
curve $C_k$, and consider two curves $C_{k+1}$,
$\tilde C_{k+1}$ as in (a). Let $P_{k+1}\in C_{k+1}$,
$Q_{k+1}\in C_{k+1}$ and the tangent line $t_{k+1}$ of $C_{k+1}$ at
$Q_{k+1}$ be as above; denote the corresponding data for $\tilde
C_{k+1}$ by $\tilde P_{k+1},$ $\tilde Q_{k+1}$ and $\tilde t_{k+1}.$
Consider the Cremona transformations $\varphi_k:=\psi(P_{k+1},
Q_{k+1},t_{k+1})$ and
$\tilde \varphi_k:=\psi(\tilde P_{k+1}, \tilde Q_{k+1},\tilde t_{k+1}),$
and also the proper transforms $C_k:=\varphi_k(C_{k+1})$ and $\tilde
C_k:=\tilde \varphi_k(\tilde C_{k+1})$. Reversing the above arguments it i=
s
easily seen that the both curves $C_k,\,\tilde C_k$ are as in (a).
By the induction hypothesis, they differ by a projective transformation
$f:\p^2\to\p^2,$ i.e.\ $f(C_k)=\tilde C_k$. For $k>0$ the points
$Q_k\in C_k$, $S_k\in C_k$ and the line  $\ell_k$ are intrinsically
defined by the curve $C_k,$ and so, $f$ maps  these data onto the
corresponding data $\tilde Q_k$, $\tilde S_k$ and $\tilde \ell_k$ for the
curve $\tilde C_k$. Moreover, in the case $k=0$ it is easily seen
that one can choose
$f$ in such a way that $f(Q_0)=\tilde Q_0$. Then again
$f(S_k)=f(\tilde S_k)$ and $f(\ell_k)=\tilde\ell_k$.  Hence, the map
$f$ is compatible with the Cremona transformations
$\varphi_k^{-1}=\psi(S_k,Q_k,\ell_k)$ and $\tilde \varphi_k^{-1}=
\psi(\tilde S_k,\tilde Q_k,\tilde \ell_k)$, i.e.\ there is a linear
transformation $g$ of $\p^2$
such that $\varphi_k\circ  g=f\circ \tilde \varphi_k$. Clearly, $g$
transforms $C_{k+1}$ into $\tilde C_{k+1}$.
\qed

\bigskip

\noindent {\bf Remarks.}  (1) By the same approach
as in the proof of Proposition 1.2,
it is possible to show the existence and uniqueness of the
rational cuspidal curves of type $(d,d-2)$ with at least three cusps,
which was done by a different method in [FlZa 2]. By the result of
loc.cit\ such a curve $C$ has exactly three cusps, say $Q,$ $P,$ $R,$
with the multiplicity sequences $(d-2),$ $(2_a),$ $(2_b),$ respectively,
where $a+b=d-2.$  Set $\ell_P:=\ol{QP},$ $\ell_R:=\ol{QR}.$ and
denote by $t_Q$ the tangent line at $Q.$  By Bezout's Theorem,
$t_Q$ intersects $C$ in one further point $S$ different from $Q.$
Performing the Cremona transformation $\psi(S,Q,\ell_P)$ to the
curve $C,$  we obtain a curve of degree
$d+1$ with the multiplicity sequences $(d-1),$ $(2_{a+1}),$ $(2_b)$ at
the cusps. Similarly, under the Cremona transformation
$\psi(P,Q,\ell_R)$  the curve $C$ is transformed into a cuspidal
curve of the same degree
$d$ with the multiplicity sequences $(d-2),$ $(2_{a+1}),$
$(2_{b-1}).$  Thus, starting from the rational cuspidal quartic
with three cusps, we can construct all such curves. It follows from
this construction that these curves are rectifiable.

\smallskip

(2) Using the above arguments, it is also possible to classify the
rational cuspidal curves of degree five with at least three cusps, which
was done by M. Namba by a different method, see \cite[Thm.2.3.10]{Nam}.

Indeed, if the largest multiplicity of a cusp is $3,$
then projecting $C$ from this point, say $Q,$ gives a two-sheeted covering
$C \to \p^1$ with two ramification points. Hence, in this case $C$ has
three
cusps, with multiplicity sequences $(3)$ (at $Q$), $(2_2),\,\,(2),$
respectively.

If all the cusps are of multiplicity $2,$ then $C$ has singular points
$P,Q,R,\ldots$ with multiplicity sequences $(2_p),$ $(2_q),$
$(2_r),\ldots,$ where $p+q+r+\ldots=6.$ We may assume that
$p\ge q\ge r\ldots.$ Projecting from $P$ gives a three-sheeted covering $C
\to \p^1$ with four ramification points. Hence, $C$ has at most
four cusps. The possibilities are as follows:

(1) $C$ has 3 cusps of type $P=(2_2),$ $Q=(2_2),$ $R=(2_2).$

(2) $C$ has 3 cusps of type $P=(2_4),$ $Q=(2),$ $R=(2).$

(3) $C$ has 3 cusps of type $P=(2_3),$ $Q=(2_2),$ $R=(2).$

(4) $C$ has 4 cusps of type $P=(2_3),$ $Q=(2),$ $R=(2),$ $S=(2).$

(5) $C$ has 4 cusps of type $P=(2_2),$ $Q=(2_2),$ $R=(2),$ $S=(2).$

\noindent
Curves as in (1) and (4) do exist and can be constructed by Cremona
transformations. The other cases are not possible, as can be seen by
the following arguments.

(5) can be excluded since
the dual curve would be a cubic with two cusps, which is impossible.

To exclude (3), denote by $t_P$ the tangent line of $C$ at $P.$ By the
Cremona transformation $\psi:= \psi(Q,P,t_P)$ a curve $C$ as in (3)
is transformed into a quartic $C'$ with three simple cusps
$P'$, $Q'$, $R'.$ It can be seen that there is a tangent line at a smooth
point $S'$ of $C'$ passing through one of the cusps, say $Q'.$ Projecting
from $Q'$
gives a two-sheeted covering $C' \to \p^1$ with three ramification points,
namely $P',$ $R'$ and $S'.$ This contradicts the Hurwitz formula.

In the case (2), consider
the blow up at $P,$ and perform an elementary transformation at
the point of the proper transform of $C$ over $P.$ Then the image of $P$
will be a point with the multiplicity sequence $(2_2).$
Performing at this point another elementary transformation and blowing
down to $\p^2,$ we arrive at the same configuration as above. Hence,
also (2) is impossible. (This last transformation may also be considered
as a Cremona transformation, namely in the points $P$, $P'$ and $P'',$
where $P'$ is infinitesimally near to $P$ and $P''$ is infinitesimally
near to $P'.$)

Similarly, using Cremona transformations for the cases 1 and 4, one
can construct these curves and show that they are rectifiable and
projectively unique. It is also possible to treat in the same way the
rational cuspidal quintics with one or two cusps.

\medskip

Finally, we give an alternative proof for the existence and uniqueness
statements of Proposition 1.2. It provides a way of computing an
explicit parameterization for these curves.
\bigskip

\noindent {\it Alternative proof of Proposition 1.2 (a)-(c). } For $k=1$
the result is known (see e.g.
\cite{Nam}).  Let $C_k\,\,(k > 1)$ be a rational cuspidal plane curve of
degree $d = 2k+3$ with three cusps $P,\,Q,\,R$ of types $(3_k),$
$(2k,\,2_k)$ and $(2),$ respectively. Since, by Bezout's Theorem, they are
not at the same line,  we may chose them as
$Q\,(0:0:1),\,P\,(0:1:0),\,R\,(1:0:0).$  We may also chose a parameterizatio=
n
$\p^1 \to C_k$ of $C_k$ such that
$(0:1) \mapsto Q,\,(1:0) \mapsto P,\,(1:1) \mapsto R.$  Then, up to
constant factors, this parameterization can be written as
$$(x,\,y,\,z) =
(s^{2k}t^3,\,\,\,\,s^{2k}(s-t)^2(as+bt),\,\,\,\,t^3(s-t)^2q_k(s,\,t))\,,$$
where $q_k \in \C [s,\,t]$ is a homogeneous polynomial of degree $2k-2.$
Let $\G$ denotes a curve parameterized as above (with $q$ instead of
$q_k$). It is enough to prove the following

\ss

\noindent {\bf Claim.} {\it There exists unique polynomials $as+bt$ and $q$
with rational coefficients, where $q(1,\,0) = 1,$
such that the multiplicity sequences of $\G$ at the points $P,\,Q,\,R \in \G=
$
start, respectively, with $(3_k),\,(2k,\,2_k)$ and $(2)$}.

\ss

Indeed, if this is the case, then, by the genus formula,
these multiplicity sequences actually coincide resp. with
$(3_k),\,(2k,\,2_k)$ and $(2),$ and so,
$C_k = \G$ up to projective equivalence. This will prove the existence of
the curves $C_k$ defined over $\Q$ for all $k > 1,$ as well as
their uniqueness, up to projective equivalence.

\ss

\noindent {\it Proof of the claim.}
It is easily seen that, after blowing up at $Q,$ the infinitesimally near
point $Q'$ to $Q$ at the proper transform $\G'$ of $\G$ will be a singular
point of multiplicity $2$ iff
$as+bt = 2s+t.$ By \cite[(1.2)]{FlZa 2}, under this condition
the multiplicity sequence of $\G$ at $Q$ starts with $(2k,\,2_k).$

In the affine chart $(\hx,\,\hz):=(x/y,\,z/y)$ centered at $P$ we have
$$\hx = {t^3 \over (s-t)^2(2s+t)},\,\,\,\,\,\,\hz={t^3q(s,\,t) \over
s^{2k}(2s+t)}\,.$$
In the sequel we denote by the same letter $t$ the affine coordinate $t/s$
in $\p^1 \setminus \{(0:1)\}.$
Thus, in this affine chart in $\p^1$ centered at $(1:0)$  we have
$$(\hx,\,\hz)= \left({t^3 \over (t-1)^2(t +2)},\,\,{t^3 \over (t
+2)}\,\hq(t)\right)\,,$$
where $\hq(t) = \sum_{i=0}^{2k-2} c_it^i$ and where, by the above
assumption, $c_0 = 1.$

After blowing up at $P,$ in the affine chart with the coordinates
$(u,\,v),$ where $(\hx,\,\hz) = (u,\,uv),$ we will have
$$(u,\,v) = (\hx,\,\hz/\hx) = \left({t^3 \over (t-1)^2(t
+2)},\,\,\,\,\hq(t)(t-1)^2\right)\,.$$  To move the origin to the
infinitesimally near point $P' \in \G'$ of $P,$
we set
$$(\hu,\, \hv) = (u,\,v-1) = \left({t^3 \over (t-1)^2(t +2)},\,\,\,\,
\hq(t)(t-1)^2 - 1\right)\,.$$
The following conditions guarantee that the multiplicity of the curve $\G'$ =
at
$P'$ is at least $3$:
$$t^3 \,\vert\,\,[\hq(t)(t-1)^2 - 1] \Longleftrightarrow$$
$$[\hq(t)(t-1)^2 - 1]'_0 = [\hq(t)(t-1)^2 - 1]''_0 =
0\Longleftrightarrow$$
$$\hq'(0) = 2,\,\hq''(0) = 6 \Longleftrightarrow c_1 = 2,
\,c_2=3\,.\eqno{(15)}$$
In the case when $k=2$ this uniquely determines the polynomial $q$:
$$q(s,\,t) = s^2 + 2st + 3t^2\,.$$  In what follows we suppose that
$k > 2.$
Assume that the conditions (15) are fulfilled. Then we have the following
coordinate presentation of $\G'$:
$$(\hu,\,\hv) = \left({t^3 \over (t-1)^2(t +2)},\,\,\,\,t^3h(t)\right)\,,$$
where $h(t):=[\hq(t)(t-1)^2 - 1] / t^3$ is a polynomial of
degree $2k-3,$ which satisfies the conditions
$$(t-1)^2 \,\vert\,\,[t^3h(t) + 1] \Longleftrightarrow
h(1) = -1,\,h'(1) = 3\,.\eqno{(15')}$$
Once ($15'$) are fulfilled, one can find $\hq$ as
$\hq = [t^3h(t) + 1]/(t-1)^2,$ and we have $\hq \in \Q[t]$ iff $h \in \Q[t=
].$

Let $\xi \in \C[[t]]$ be such that
$\xi^3 = {t^3 \over (t-1)^2(t +2)}.$ By \cite[(3.4)]{FlZa 2}, the
multiplicity sequence of $\G'$ at $P'$ starts with $(3)_{k-1}$ iff
$$t^3h(t) \equiv \hf(\xi^3) \,\,\,{\rm mod}\,\xi^{3(k-1)}\,,$$
where $\hf = \sum\limits_{i=0}^{k-1} \ha_ix^i \in \C[x]$ is a polynomial
of degree $ \le k-1.$
Multiplying the both sides by the unit $[(t-1)^2(t +2)]^{k-1} \in \C[[t]]$,
we will get
$$[(t-1)^2(t +2)]^{k-1}t^3h(t) \equiv [(t-1)^2(t +2)]^{k-1}
\sum\limits_{i=0}^{k-1} \ha_i\xi^{3i} \equiv \sum\limits_{i=0}^{k-1}
\ha_it^{3i}[(t-1)^2(t +2)]^{k-1-i} \,\,\,\,{\rm mod}\,t^{3(k-1)}\,.$$
Since, by our assumption, $k > 1,$ we should have $\ha_0 = 0,$ and after
dividing out the factor $t^3,$ we get
$$[(t-1)^2(t +2)]^{k-1}h(t) \equiv \sum\limits_{i=0}^{k-2}
\ha'_it^{3i}[(t-1)^2(t +2)]^{k-2-i} \,\,\,\,{\rm mod}\,t^{3(k-2)}\,,$$
where $\ha'_i = \ha_{i-1},\,\,i=1,\dots,k-2.$  In other words, we have
$$[(t-1)^2(t +2)]^{k-1}h(t) = \hf(t^{3}, \,(t-1)^2(t +2))
+ \hg(t)t^{3(k-2)}\,,$$
where $\hf(x,\,y) = \hf_k(x,\,y) := \sum_{i=0}^{k-2} \ha'_ix^iy^{k-2-i=
}$ is
a homogeneous polynomial of degree $k-2,$ and hence $\hg(t) = \hg_k(t) =
\sum_{i=0}^{2k} \hb_it^i$
should be a polynomial of degree $2k.$
Denoting $\tau = t^3$ and
$\lambda = (t-1)^2(t +2) = t^3 - 3t + 2$,
we have
$$\lambda^{k-1}h = \hf(\tau,\,\lambda) + \tau^{k-2}\hg\,.$$  Observe that
$\hf(\tau,\,\lambda)$ (resp. $\tau^{k-2}\hg$)
contains the monomial $\ha'_0\tau^{k+2}$ (resp. $\hb_0\tau^{k+2}$).
To avoid indeterminacy, we may assume, for instance, that $\ha'_0=0.$
Then $\hf = \lambda f,$ where $f(x,\,y) := \sum_{i=0}^{k-2}
a_ix^iy^{k-3-i},\,\,\, a_i := \ha'_{i-1},\,i=0,\dots, k-3,$ and so
$$\lambda^{k-1}h = \lambda f(\tau,\,\lambda) + \tau^{k-2}\hg\,.$$
Since $(\tau,\,\lambda) = 1,$ we have $\lambda \,\vert \,\hg,$ that is,
$\hg = \lambda g,$ where $g(t) := \sum_{i=0}^{2k-3} b_it^i.$
Finally, we arrive at the relation
$$\lambda^{k-2}h(t) = f(\tau,\,\lambda) + \tau^{k-2}g(t)\,,$$
where deg$\,f = k-3,$ deg$\,h= $deg$\,g = 2k-3,$ and $h$ should satisfy
the conditions ($15'$). It follows that
$$\lambda^{k-2}\,\vert\,\,[f(\tau,\,\lambda) + \tau^{k-2}g]\,,\eqno{(16)}$$
and $$\tau^{k-2}\,|\,\,[f(\tau,\,\lambda) -\lambda^{k-2}h]\,.\eqno{(16')}$$
Each of these conditions together with ($15'$) determines the triple of
polynomials $f,\,g,\,h$ as above in a unique way. Indeed, once
$f$ and $g$ satisfy ($15'$) and (16), we can find $h$ as
$h = [f(\tau,\,\lambda) + \tau^{k-2}g]/ \lambda^{k-2}.$ Actually, (16) is
equivalent to the vanishing of derivatives of the function
$f(\tau,\,\lambda) + \tau^{k-2}g \in \C[t]$
at the point $t = 1$ up to order $2k-5$ and at the point $t = -2$ up to
order $k-3.$ This yields a system of $3k-6$ linear equations in the $3k-4$
unknown coefficients of $f$ and $g$; ($15'$) provides another
two linear equations. That is, we have the following system:
$$
\left(f(\tau,\,\lambda) + \tau^{k-2}g\right)^{(m)}_{t=-2} =
0,\,\,\,m=0,\dots,k-3$$
$$\left(f(\tau,\,\lambda) + \tau^{k-2}g\right)^{(m)}_{t=1} =
0,\,\,\,m=0,\dots,2k-5$$
$$\left(f(\tau,\,\lambda) + \tau^{k-2}g\right)^{(2k-4)}_{t=1} =
-3^{k-2}(2k-2)!\eqno{(S)}$$
$$\left(f(\tau,\,\lambda) + \tau^{k-2}g\right)^{(2k-3)}_{t=1} =
-3^{k-3}(k-11)(2k-1)!$$
(Indeed, put $u = t-1$; in view of ($15'$) we have
$$\lambda = (t-1)^2(t+2) = u^2(u+3),\,\,\,\,h(t) = -1 + 3u + \dots\,,$$
and hence
$$f(\tau,\,\lambda) + \tau^{k-2}g(t) = \lambda^{k-2}h(t) =
[u^2(u+3)]^{k-2}h(t) =$$
$$ u^{2k-4}(3^{k-2} + (k-2)3^{k-3}u + \dots)(-1 + 3u + \dots) =
u^{2k-4}(-3^{k-2} - 3^{k-3}(k-11)u + \dots)\,.)$$
The system ($S$) has a unique solution iff it is so for the associated
homogeneous system, say, ($S_0$). Passing from ($S$) to ($S_0$) actually
corresponds to passing from $h$ to a polynomial
$h_0$ of degree $\le 2k-3$ which satisfies, instead of ($15'$),
the conditions
$$h_0(1) = h'_0(1) = 0 \Longleftrightarrow (t-1)^2 \,\vert\,\,h_0(t)
\Longleftrightarrow h_0(t) =
(t-1)^2 {\widetilde h}(t),\,\,\,{\rm deg}\,{\widetilde h} \le 2k-5\,.
\,\eqno{(15'')}\,.$$  Thus, we have to prove that the equality
$$\lambda^{k-2}(t-1)^2 {\widetilde h}(t) = f(\tau,\,\lambda) +
\tau^{k-2}g(t)\,,$$
where $f=0$ or deg$\,f = k-3,$ deg$\,g \le 2k-3,$ and
${\rm deg}\,{\widetilde h} \le 2k-5,$
is only possible for $f=g={\widetilde h}=0.$ Or, equivalently, we have=
 to
show that
the $5k-8$ polynomials in $t$ in the union $T$ of the three systems:
$$T_1:= \left\{\tau^i\lambda^{k-3-i}\right\}_{i=0,\dots,k-3},\,\,\,
T_2:=\left\{t^i(t-1)^2\lambda^{k-2}\right\}_{i=0,\dots,2k-5},\,\,\,
T_3:=\left\{t^i\tau^{k-2}\right\}_{i=0,\dots,2k-3}$$
are linearly independent. After replacing the system $T_2$ by the
equivalent one:
$$T'_2:= \left\{(t-1)^{2k-2}(t+2)^{k-2+i}\right\}_{i=0,\dots,2k-5}\,,$$
we will present these three systems as follows:
$$T_1=\left\{p_i:=\tau^{k-3-i}\lambda^i=t^{3(k-2-i)}(t-1)^{2i}(t+2)^i,
\,\,\,\,i=0,\dots,k-3\right\}$$
$$T'_2=\left\{p_i:=(t-1)^{2k-2}(t+2)^i,
\,\,\,\,i=k-2,\dots,3k-7\right\}$$
$$T_3=\left\{p_i:=t^i,\,\,\,\,i=3k-6,\dots,5k-9\right\}\,.$$  Denote
 $P= $
span$\,(T_1,\,T_2,\,T_3)= $ span$\,(T_1,\,T'_2,\,T_3).$
Note that deg$\,p \le 5k-9$ for all $p \in P,$ that is, dim$\,P \le 5k-8.$
Consider the following system of $5k-8$ linear functionals on $P$:
$$\varphi_i\,:\,p \longmapsto p^{(i)}(-2),\,\,\,\,i=0,\dots,3k-7\,,$$
$$\varphi_i\,:\,p \longmapsto p^{(i)}(0),\,\,\,\,i=3k-6,\dots,5k-9\,.$$
It
is easily seen that the matrix
$M:=\left(\varphi_i(p_j)\right)_{i,\,j=0,\dots,5k-9}$
is triangular with non-zero diagonal entries. This proves that, indeed,
rang$\,T =$dim$\,P=5k-8,$ as stated.

The coefficients of the system ($S$) being integers, its unique solution is
rational,
i.e.\ the polynomials $f$ and $g$ are defined over $\Q.$  It follows as
above that the polynomials $h$ and $q$ are also defined over
$\Q.$
This completes the alternative proof of Proposition 1.2. \qed

\bigskip

\noindent {\bf Remarks.} (1) In principle, the method used in the proof
allows to compute
explicitly parameterizations of the curves $C_k.$  For instance,
we saw above that
for $k=2$ a parameterization of $C_2$ is given by the choice
$$
q_2(s,t):=s^2+2st+3t^2,\quad a:=2,\quad b:=1\,.
$$

\noindent (2) We have to apologize for a pity mistake in Lemma 4.1(b)
[FlZa 2, Miscellaneous] (this does not affect the other results of [FlZa 2],
besides only the immediate Corollary 4.2).

\bigskip\footnotesize
\renewcommand{\arraystretch}{0.8}
\hspace{2mm}\begin{tabular}{lccccl}

 Hubert Flenner    &&&&& Mikhail Zaidenberg\\
 Fakult{\"a}t f{\"u}r Mathematik  &&&&& Universit{\'e}
Grenoble I \\
Ruhr Universit{\"a}t Bochum  &&&&& Institut Fourier \\
Geb.\ NA 2/72     &&&&&    UMR 5582 CNRS-UJF \\
Universit{\"a}tsstr.\ 150    &&&&&BP 74\\
44780 BOCHUM &&&&&    38402 St. Martin
d'H{\`e}res--c{\'e}dex \\
Germany &&&&& France\\
e-mail:&&&&& e-mail:\\
Hubert.Flenner@rz.ruhr-uni-bochum.de&&&&&zaidenbe@ujf-grenoble.fr

\end{tabular}

\end{document}